\documentclass[12pt]{article} 
\def\be{\begin{equation}}
\def\ee{\end{equation}}
\def\bea{\begin{eqnarray}}
\def\eea{\end{eqnarray}}
\def\nn{\nonumber}
\def\ba{\begin{array}}
\def\ea{\end{array}}
\def\lb{\left(}
\def\rb{\right)}
\def\ls{\left[}
\def\rs{\right]}
\def\lc{\left\{}
\def\rc{\right\}}
\def\ul{\underline}
\def\s{\sum_{n=0}^\infty}
\def\ds{\sum_{n=-\infty}^\infty} 
\def\f{\frac}
\def\a{\alpha}
\def\b{\beta}
\def\D{\Delta}
\def\l{\lambda}
\def\sg{\sigma}
\def\z{\zeta}
\def\t{\theta}
\def\lra{\longrightarrow}
\begin{document}
\begin{center}
{\large\bf Two-parameter quantum algebras, twin-basic numbers, \\ 
and associated generalized hypergeometric series}\footnote{To appear in 
the Proceedings of the International Conference on Number Theory and 
Mathematical Physics, 20-21 December 2005, Srinivasa Ramanujan Centre, 
Kumbakonam, India}  

\bigskip 

R. Jagannathan\footnote{E-mail: jagan@imsc.res.in} \\
{\it The Institute of Mathematical Sciences \\ 
C.I.T. Campus, Tharamani, Chennai - 600113, India} \\
and \\
K. Srinivasa Rao\footnote{E-mail: rao@imsc.res.in} \\
{\it Srinivasa Ramanujan Centre \\ 
Shanmugha Arts, Science, Technology \& Research Academy (SASTRA) \\ 
Kumbakonam - 612001, India}  
\end{center}

\vspace{1cm}

\begin{abstract}
We give a method to embed the $q$-series in a $(p,q)$-series 
and derive the corresponding $(p,q)$-extensions of the known $q$-identities.  
The $(p,q)$-hypergeometric series, or twin-basic hypergeometric series 
(different from the usual bibasic hypergeometric series), is based on 
the concept of twin-basic number $[n]_{p,q}$ $=$ $(p^n-q^n)/(p-q)$.  
This twin-basic number occurs in the theory of two-parameter quantum 
algebras and has also been introduced independently in combinatorics. 
The $(p,q)$-identities thus derived, with doubling of the number of 
parameters, offer more choices for manipulations; for example, results 
that can be obtained via the limiting process of confluence in the usual 
$q$-series framework can be obtained by simpler substitutions.  The 
$q$-results are of course special cases of the $(p,q)$-results corresponding 
to choosing $p$ $=$ $1$.  This also provides a new look for the 
$q$-identities. 
\end{abstract}

\vspace{1cm}

\noindent
{\bf 1. Introduction} \\ 

\noindent
For the two-parameter quantum group $GL_{p,q}(2)$ the fundamental 
representation is given by the $T$-matrix,  
\be
T = \lb\ba{cc}
        a & b \\
        c & d 
    \ea\rb,
\ee
whose elements satisfy the commutation relations
\bea
   &   & ab = p^{-1}ba, \quad cd = p^{-1}dc, 
   \quad ac = q^{-1}ca, \quad bd = q^{-1}db, \nn \\ 
   &   & bc = q^{-1}pcb, \quad ad-da = (p^{-1}-q)bc, 
\eea 
consistent with the equation 
\be
R(T\otimes I)(I\otimes T) = (I\otimes T)(T\otimes I)R, 
\ee
corresponding to the $R$-matrix
\be
R = (pq)^{1/4}\lb\ba{cccc}
    (pq)^{-1/2} &      0                 &     0       & 0           \\
          0     & (p/q)^{-1/2}           &     0       & 0           \\
          0     & (pq)^{-1/2}-(pq)^{1/2} & (p/q)^{1/2} & 0           \\
          0     &      0                 &     0       & (pq)^{-1/2} 
    \ea\rb.
\ee 
The two-parameter quantum algebra, $U_{p,q}(gl(2))$, dual to $GL_{p,q}(2)$,  
is generated by $\{Z, J_0, J_\pm\}$ satisfying the commutation relations 
\bea
   &   & [Z,J_0] = 0, \quad [Z,J_\pm] = 0, \nn \\
   &   & [J_0,J_\pm] =  \pm J_\pm, \quad 
J_+J_--pq^{-1}J_-J_+ = \f{p^{-2J_0}-q^{2J_0}}{p^{-1}-q}.
\label{pqgl2}
\eea 
To realize this algebra~(\ref{pqgl2}), a $(p,q)$-oscillator algebra, 
\be
aa^\dagger-qa^\dagger a = p^{-N}, \quad 
[N,a] = -a, \quad [n,a^\dagger] = a^\dagger, 
\label{pqosc}
\ee
was introduced in~\cite{CJ} generalizing/unifying several forms of 
$q$-oscillator algebras well known in the earlier physics literature 
related to the representation theory of single-parameter quantum 
algebras.  The algebra~(\ref{pqosc}) is satisfied when 
\be
a^\dagger a = \f{p^{-N}-q^N}{p^{-1}-q}, \qquad 
aa^\dagger = \f{p^{-(N+1)}-q^{N+1}}{p^{-1}-q}.
\label{aaN}
\ee  
When $p$ $=$ $q$ or $p$ $=$ $1$ the algebra~(\ref{pqosc}) becomes two 
different versions of the $q$-oscillator algebra related to the 
representation theory of $U_q(sl(2))$. 

The relations~(\ref{pqgl2} and~(\ref{aaN}) suggest immediately a 
generalization of the Heine $q$-number, 
\be
[n]_q = \f{1-q^n}{1-q}, 
\label{nq}
\ee 
to a $(p,q)$-number as 
\be
[n]_{p,q} = \f{p^n-q^n}{p-q}.  
\label{npq}
\ee
If we define a $(p,q)$-derivative by 
\be
\hat{D}_{p,q}f(z) = \f{f(pz)-f(qz)}{(p-q)z}  
\label{2pdd}
\ee 
then 
\be
\hat{D}_{p,q}z^n = [n]_{p,q}z^{n-1}.
\ee
Several properties of this $(p,q)$-number~(\ref{npq}), which we will now  
call as the twin-basic number, including the elements of $(p,q)$-calculus 
following from~(\ref{2pdd}) were studied very briefly in~\cite{CJ}.  For 
the sake of convenience, we shall denote $[n]_{p,q}$ simply as $[n]$ and 
omit the subscripts $p,q$ from other expressions also whenever the values 
of these twin-base parameters are clear from the context.  

Around the same time as~\cite{CJ}, Brodimas, {\em et al}.~\cite{Bet} and 
Arik, {\em et al}.~\cite{Aet} also, independently, introduced the 
$(p,q)$-number in the physics literature, but in a very much less detailed 
manner.  They also introduced the $(p,q)$-oscillator and the $(p,q)$-number 
in the same context of realization of $U_{p,q}(gl(2))$.  It is a surprising 
fact that around the same time, without any connection to the quantum group 
related mathematics/physics literature, Wachs and White~\cite{WW} introduced 
the $(p,q)$-number, defined as $(p^n-q^n)/(p-q)$, in the mathematics 
literature while generalizing the Sterling numbers, motivated by certain 
combinatorial problems (for further generalizations and applications in 
this direction see~\cite{RW}).  In physics literature, Katriel and 
Kibler~\cite{KK} defined the $(p,q)$-binomial coefficients and derived a 
$(p,q)$-binomial theorem while discussing normal ordering for deformed boson 
operators obeying the algebra~(\ref{pqosc}).  Smirnov and Wehrhahn~\cite{SW} 
gave an operator, or noncommutative, version of such a $(p,q)$-binomial 
theorem.  Floreanini, Lapointe and Vinet~\cite{FLV} related the 
algebra~(\ref{pqosc}) to bibasic hypergeometric functions~\cite{AV1,AV2}.  
Burban and Klimyk~\cite{BK} studied the $(p,q)$-differentiation, 
$(p,q)$-integration, and the $(p,q)$-hypergeometric series $_r\Psi_{r-1}$ 
in detail.  Gelfand, {\em et al}.~\cite{Get1,Get2} generalized the 
two-parameter deformed derivative~(\ref{2pdd}) and developed a very general 
theory of deformation of classical hypergeometric functions.  Their general 
formalism of deformed hypergeometric functions is based on a $u$-derivative 
\be
\hat{D}_uf(z) = \f{1}{z}u\lb z\f{d}{dz}\rb f(z)
\ee 
where $u(z)$ is an arbitrary entire function.  This leads to a $u$-calculus 
and a unified exposition of the classical theory and the $q$-theory and 
results in new $u$-analogues of classical hypergeometric functions.  The 
$(p,q)$-hypergeometric series corresponds to the choice $u(z)$ $=$ 
$(p^z-q^z)/(p-q)$.  Generalizing the definition of $_r\Psi_{r-1}$ by 
Burban and Klimyk~\cite{BK}, one of us defined the general 
$(p,q)$-hypergeometric series $_r\Phi_s$ and derived some related 
preliminary results~\cite{J}. Some applications of the 
$(p,q)$-hypergeometric series in the context of representations of 
two-parameter quantum groups have been considered by Nishizawa~\cite{N} 
and Sahai and Srivastava~\cite{SS}.  

In the present work we shall deal only with the $(p,q)$-hypergeometric 
series as defined in~\cite{J}.  We introduce a method of application of 
the $(p,q)$-series to convert the various well known $q$-identities into 
their $(p,q)$-analogues; after the conversion the resulting 
$(p,q)$-identities offer more choices for symbolic manipulations 
transcending the applications of the original $q$-identities and in fact 
give a new look to the latter.  \\

\noindent
{\bf 2. Twin-basic hypergeometric series $_r\Phi_s$} \\

\noindent 
Let us recall some basic definitions from the theory of $q$-hypergeometric 
series~\cite{GR}.  The $q$-shifted factorial is given by 
\bea
(a;q)_n & = & \lc\ba{l}
                  1, \quad n = 0, \\ 
        (1-a)(1-aq)(1-aq^2)\dots(1-aq^{n-1}), \\ 
        \quad \quad n = 1,2,\dots.  
\ea\right. 
\label{aqn}
\eea
With 
\be
(a_1,a_2,\dots,a_k;q)_n = (a_1;q)_n(a_2;q)_n\dots(a_k;q)_n. 
\label{aan}
\ee
the $q$-hypergeometric series, or the basic hypergeometric series,  
$_r\phi_s$ is defined as 
\bea 
   &   & _r\phi_s(a_1,a_2,\ldots,a_r;b_1,b_2,\ldots,b_s;q,z) \nn \\
   &   &  \quad = \s\f{(a_1,a_2,\ldots,a_r;q)_n}          
                              {(b_1,b_2,\ldots,b_s;q)_n(q;q)_n}
                              ((-1)^nq^{n(n-1)/2})^{1+s-r}z^n. 
\label{phi}
\eea
Let us now call the $(p,q)$-number~(\ref{npq}) as twin-basic number and 
define the twin-basic analogues of~(\ref{aqn}) and~(\ref{aan}) as follows: 
\bea
((a,b);(p,q))_n & = & \lc 
\ba{l}
1, \qquad n = 0, \\
(a-b)(ap-bq)(ap^2-bq^2)\dots(ap^{n-1}-bq^{n-1}),   \\ 
   \quad \qquad n = 1,2,\dots. 
\ea\right. 
\eea 
\bea
   &   & ((a_{1p},a_{1q}),( a_{2p},a_{2q}),\dots,
(a_{mp},a_{mq});(p,q))_n \nn \\
   &   &  \quad = ((a_{1p},a_{1q});(p,q))_n 
((a_{2p},a_{2q});(p,q))_n \dots ((a_{mp},a_{mq});(p,q))_n. 
\eea
Note that 
\be
(a;q)_n = ((1,a);(1,q))_n.
\ee 
Then, the $(p,q)$-analogue of~(\ref{phi}), the $(p,q)$-hypergeometric 
series, or the twin-basic hypergeometric series, can be defined as  
\bea
   &   & 
_r{\Phi}_s((a_{1p},a_{1q}),\dots,
(a_{rp},a_{rq});(b_{1p},b_{1q}),\dots, 
(b_{sp},b_{sq});(p,q),z) \nn \\ 
   &   & = \s\,\f{((a_{1p},a_{1q}),\dots,
(a_{rp},a_{rq});(p,q))_n}{((b_{1p},b_{1q}),\dots,
(b_{sp},b_{sq});(p,q))_n ((p,q);(p,q))_n} \nn \\
   &   & \qquad \qquad \qquad \qquad \qquad 
\times((-1)^n(q/p)^{n(n-1)/2})^{1+s-r}z^n,
\label{Phi}
\eea
with $|q/p|<1$~\cite{J}.  Though, generally, we shall assume $0 < q < p$, 
$p$ and $q$ can also take other values if there is no problem with 
convergence of the particular series involved in a result.  When 
$a_{1p}=a_{2p}=\dots=a_{rp}$$=b_{1p}=b_{2p}=\dots=b_{sp}=1$, 
$a_{1q}=a_1,a_{2q}=a_2,\dots,a_{rq}=a_r$, $b_{1q}=b_1,b_{2q}=b_2,\dots, 
b_{s,q}=b_s$, and $p=1$, $_r\Phi_s$ $\lra$ $_r\phi_s$.  Special interesting 
choices for $(p,q)$, from the point of view of quantum groups, are 
$(q^{-1/2},q^{1/2})$, $(q^{-1},q)$ and, more generally, $(p^{-1},q)$.  
Throughout the paper we shall assume $|z|<1$.  Also, we shall assume 
all the parameters to be generic, with nonzero values, unless specified 
otherwise.  While referring to the classical results of the $q$-series 
we shall use the standard notations as in~\cite{GR} (see also~\cite{K}).  
Often, the parameter doublets $(a_p,a_q)$, $(b_p,b_q)$, etc., will be 
denoted by different symbols according to the convenience of the situation 
and such notations should be clear from the context. 

Let us recall the definition of a bibasic hypergeometric series with two 
bases $q$ and $q_1$~\cite{AV1,AV2} (see also~\cite{GR}): 
\bea
  &   &  {\cal F}(\ul{a},\ul{c};\ul{b},\ul{d};q,q_1,z) = \nn \\ 
  &   & \qquad  
\s\f{(\ul{a};q)_n(\ul{c};q_1)_n}{(\ul{b};q)_n(\ul{d};q_1)_n(q;q)_n} \nn \\ 
  &   &  \qquad \times 
         ((-1)^nq^{n(n-1)/2})^{1+s-r}((-1)^nq_1^{n(n-1)/2})^{s_1-r_1}z^n, 
\label{bibase}
\eea
where $\ul{a}$ $=$ $(a_1,a_2,\dots,a_r)$, 
$\ul{b}$ $=$ $(b_1,b_2,\dots,b_s)$, $\ul{c}$ $=$ $(c_1,c_2,\dots,c_{r_1})$,
and $\ul{d}$ $=$ $(d_1,d_2,\dots,d_{s_1})$.  It is clear that 
in~(\ref{bibase}) the two unconnected bases $q$ and $q_1$ are regarded 
are assigned partially to different numerator and denominator parameters 
whereas in the twin-basic hypergeometric series~(\ref{Phi}) the twin 
base parameters $p$ and $q$ are inseparable and assigned to all the numerator 
and denominator parameter doublets.  

Let 
\be
\D_{(\a,\b)}f(z) = \a f(qz) - \b f(pz). 
\ee
With 
\be
\D f(z) = \D_{(1,1)}f(z) = f(qz) - f(pz),  
\ee 
it may be noted that 
\be
\hat{D}f(z) = \f{\D f(z)}{\D z}.
\ee
Then it is seen that $_r\Phi_s$ satisfies the $(p,q)$-difference equation 
\be
\lb \D\prod_{i=1}^s\D_{\lb b_{iq}/q,b_{ip}/p\rb}\rb{}_r\Phi_s = 
\lb z\prod_{i=1}^r\D_{\lb a_{iq},a_{ip}\rb}\rb
{}_r\Phi_s\lb (q/p)^{1+s-r}z\rb.
\ee
When $a_{1p} = a_{2p} =\dots= a_{rp}$ $=b_{1p} = b_{2p} =\dots = 
b_{sp} = 1$, $a_{1q} = a_1, a_{2q} = a_2,\dots, a_{rq} = a_r$, 
$b_{1q} = b_1, b_{2q} = b_2,\dots,b_{s,q} = b_s$, and $p = 1$ this equation 
reduces to the $q$-difference equation satisfied by $_r\phi_s$.  

Let us now construct a method to embed the usual 
$_r\phi_s$-series~(\ref{phi}) in the $_r\Phi_s$-series~(\ref{Phi}).  To 
this end, we note 
\be
((\l a,\l b);(p,q))_n = \l^n((a,b);(p,q))_n,
\label{ll}
\ee
for any arbitrary nonzero $\l$, and 
\be
(b/a;q/p)_n = a^{-n}p^{-n(n-1)/2}((a,b);(p,q))_n. 
\label{sf12}
\ee 
Thus, we can write, formally, 
\bea
   &   & _r\phi_s(a_{1q}/a_{1p},a_{2q}/a_{2p},..,a_{rq}/a_{rp};
b_{1q}/b_{1p},b_{2q}/b_{2p},..,b_{sq}/b_{sp};q/p,z) \nn \\
   &   & = \lc 
\ba{l}
_r\Phi_s((a_{1p},a_{1q}),..,(a_{rp},a_{rq});
               (b_{1p},b_{1q}),..,(b_{sp},b_{sq});(p,q),\mu z)\\
\qquad \qquad \qquad \qquad \qquad \quad\ \ {\rm if}\ s = r-1, \\
_{s+1}\Phi_s((a_{1p},a_{1q}),..,(a_{rp},a_{rq}),(0,1),..,(0,1); 
                 (b_{1p},b_{1q}),..,(b_{sp},b_{sq});  \\ 
\qquad \qquad \qquad (p,q),\mu z), \quad {\rm if}\ s>r-1, \\ 
_r\Phi_{r-1}((a_{1p},a_{1q}),..,(a_{rp},a_{rq});(b_{1p},b_{1q}),
                ..,(b_{sp},b_{sq}),(0,1),..,(0,1);  \\
\qquad \qquad \qquad (p,q),\mu z), \quad {\rm if}\ s<r-1, 
\ea \right. \nn \\ 
   &   & \qquad \qquad \qquad \qquad {\rm with}\ 
\mu = \f{b_{1p}b_{2p}.. b_{sp}p}{a_{1p}a_{2p}..a_{rp}}, 
\label{phi12}
\eea
assuming that the given $_r\phi_s$-series is convergent or terminating. 
Hence any well behaved $\phi$-series can be written as a $\Phi$-series.  
But, the converse is not true, in general; in the general case, when 
$p\neq 1$, this is possible only for an $_r\Phi_{r-1}$.  To see this,  
it is enough to look at $_0\Phi_0$: 
\bea
_0\Phi_0(-;-;(p,q),z) & = & 
\s\f{(-1)^n(q/p)^{n(n-1)/2}}{((p,q);(p,q))_n}z^n \nn \\ 
   & = & \s\f{(-1)^n(\rho/p)^{n(n-1)/2}}{(\rho;\rho)_n}(z/p)^n, 
\quad {\rm with}\ \rho = q/p, \nn \\
   &   &  
\eea 
which shows that $_0\Phi_0$ becomes a $\phi$-series if and only if $p$ $=$ 
$1$.  Similarly, one is easily convinced that a generic $_r\Phi_s$-series 
cannot be identified within the class of $\phi$-series unless $p$ $=$ $1$ 
or $s = r-1$ (the first case in the above equation~(\ref{phi12})).  It is 
thus clear that the $(p,q)$-series is a larger structure in which the 
$q$-series gets embedded.  Also, note that in the usual theory of 
$\phi$-series there is no direct analogue for the choice $a_{ip} = 0$ or 
$b_{ip} = 0$, for any $i$, permissible, in general (of course, subject to 
conditions of convergence and so on), in the case of the $(p,q)$-series; 
to obtain a corresponding result in the case of the $\phi$-series one will 
have to resort to the limit process of confluence, namely, replacing $z$ 
by $z/a_r$ and taking the limit $a_r$ $\lra$ $\infty$.  As an example 
consider the following.  As is well known, in the definition of the usual 
$q$-hypergeometric series~(\ref{phi}), presence of the factor 
$((-1)^nq^{n(n-1)/2})^{1+s-r}$ (absent in the earlier 
literature~\cite{B,A,S}) leads to the useful relation 
\be
\lim_{a_r\lra\infty}{}_r\phi_s(z/a_r) = {}_{r-1}\phi_s(z).
\ee  
For the $(p,q)$-hypergeometric series~(\ref{Phi}) the corresponding property 
is:  
\bea 
   &   &  \lim_{a_{rq}\lra\infty}{}_r\Phi_s(z/a_{rq}) = 
{}_r{\Phi}_s((a_{1p},a_{1q}),..,(a_{(r-1)p},a_{(r-1)q})
(0,1); \nn \\ 
   &   &  \qquad \qquad \qquad \qquad 
          (b_{1p},b_{1q}),..,(b_{sp},b_{sq});(p,q),z) \nn \\ 
   &   &  \lim_{a_{rp}\lra\infty}{}_r\Phi_s(z/a_{rp}) =  
{}_r{\Phi}_s((a_{1p},a_{1q}),..,(a_{(r-1)p},a_{(r-1)q})
(1,0); \nn \\ 
   &   &  \qquad \qquad \qquad \qquad  
          (b_{1p},b_{1q}),..,(b_{sp},b_{sq});(p,q),z).
\eea              
Let us also note down the converse of~(\ref{phi12}) in the case $s = r-1$: 
\bea 
   &   & _r\Phi_{r-1}((a_{1p},a_{1q}),..,(a_{rp},a_{rq});  
         (b_{1p},b_{1q}),..,(b_{r-1,p},b_{r-1,q});(p,q),z)\nn \\
   &   & = {}_r\phi_{r-1}(a_{1q}/a_{1p},a_{2q}/a_{2p},..,a_{rq}/a_{rp}; \nn \\ 
   &   & \qquad b_{1q}/b_{1p},b_{2q}/b_{2p},..,
                 b_{r-1,q}/b_{r-1,p}; q/p,z/\mu).  
\label{phi21}
\eea
Another set of relations often useful are 
\bea
\f{(b/a;q/p)_\infty}{(d/c;q/p)_\infty} & = & 
\lim_{N\lra\infty}\f{(b/a;q/p)_N}{(d/c;q/p)_N} \nn \\ 
   & = & \lim_{N\lra\infty}\f{a^{-N}p^{-N(N-1)/2}((a,b);(p,q))_N}
{c^{-N}p^{-N(N-1)/2}((c,d);(p,q))_N} \nn \\
   & = & \f{((c,bc/a);(p,q))_\infty}{((c,d);(p,q))_\infty} \nn \\ 
   & = & \f{((a,b);(p,q))_\infty}{((a,ad/c);(p,q))_\infty}\,,
\label{lim}
\eea
and its obvious generalizations containing several factors in the 
numerator and denominator. 

Manipulations using the above relations take the usual $q$-identities to 
$(p,q)$-identities.  The original $q$-identities are, of course, special 
cases corresponding to the choice 
$a_{1p} = a_{2p} = \ldots a_{rp} = b_{1p} = b_{2p} = \ldots = b_{r-1,p} = 1$, 
and $p$ $=$ $1$.  We shall consider a few examples below. \\

\noindent
{\bf 3. $(p,q)$-Binomial theorem}  \\

\noindent
The usual $q$-binomial theorem is 
\be
_1\phi_0(a;-;q,z) = \f{(az;q)_\infty}{(z;q)_\infty}\,.
\label{qbinom}
\ee
The $(p,q)$-analogue of this is given by 
\be
_1\Phi_0((a,b);-;(p,q),z) = 
\f{((p,bz);(p,q))_\infty}{((p,az);(p,q))_\infty}. 
\label{pqbinom}
\ee
{\em Proof}~:  
Let us rewrite~(\ref{qbinom}) as 
\be
_1\phi_0(b/a;-;q/p,\z) = \f{(b\z/a;q/p)_\infty}{(\z;q/p)_\infty}.
\ee
Using~(\ref{phi12}) and~(\ref{lim}), we have 
\be
_1{\Phi}_0((a,b);-;(p,q),p\z/a) = 
\f{((a,b\z);(p,q))_\infty}{((a,a\z);(p,q))_\infty}\,.
\ee
Now, taking $\z = za/p$, we get
\be
_1{\Phi}_0((a,b);-;(p,q),z) = 
\f{((a,abz/p);(p,q))_\infty}{((a,a^2z/p);(p,q))_\infty}. 
\ee
Using the arguments of~(\ref{ll}) and~(\ref{lim}), by pulling out 
powers of $a/p$ in the numerator and denominator of the r.h.s., the 
$(p,q)$-binomial theorem~(\ref{pqbinom}) follows. 

The usual $q$-binomial theorem~(\ref{qbinom}) is recovered when $a$ $=$ $1$ 
and $p$ $=$ $1$.  The $(p,q)$-binomial theorem obtained in~\cite{BK} is a 
special case of~(\ref{pqbinom}) corresponding to the specific choice 
$(a,b)$ $=$ $(q^{-a/2},p^{a/2})$ and $(p,q)$ $=$ $(q^{-1/2},p^{1/2})$.  
An interesting feature of the $(p,q)$-binomial theorem~(\ref{pqbinom})
may be noted here. The product 
$\prod_{i=1}^n {}_1\Phi_0((a_{ip},a_{iq});-;(p,q),z)$ is seen to be an 
invariant under the group of independent permutations of the $p$-components 
$(a_{1p},a_{2p},\dots ,a_{np})$ and the $q$-components 
$(a_{1q},a_{2q},\dots ,a_{nq})$.  This product has value $1$ if the 
$n$-tuple of $p$-components $(a_{1p},a_{2p},\dots,a_{np})$ is 
related to the $n$-tuple of $q$-components $(a_{1q},a_{2q},\dots,a_{nq})$ 
by a mere permutation.  

For $n$ $=$ $2$ this result implies that 
\be
_1\Phi_0((a,b);-;(p,q),z){}_1\Phi_0((b,a);-;(p,q),z) = 1.
\label{abba}
\ee
A special case of this relation is 
\be
_1\Phi_0((1,0);-;(1,q),z){}_1\Phi_0((0,1);-; (1,q),z) = 1. 
\label{phiprod}
\ee
Recognizing that 
\bea
_1\Phi_0((1,0);-;(1,q),z) & = & \s\f{1}{(q;q)_n}z^n = e_q(z), \nn \\ 
_1\Phi_0((0,1);-;(1,q),z) & = & \s\f{q^{n(n-1)/2}}{(q;q)_n}(-z)^n 
                            = E_q(-z),   
\eea 
where $e_q(z)$ and $E_q(z)$ are the canonical $q$-exponentials, the well 
known relation 
\be
e_q(z)E_q(-z) = 1, 
\ee
follows from~(\ref{phiprod}).  It should be noted that, while in the 
usual $q$-theory~\cite{GR} $e_q(z)$ is $_1\phi_0(0;-;q,z)$ and $E_q(z)$ 
is $_0\phi_0(-;-;q,-z)$, in the $(p,q)$-series formalism both $e_q(z)$ and 
$E_q(z)$ belong to the same $_1\Phi_0$-series.  This result suggests the 
natural definitions 
\bea
e_{p,q}(z) = {}_1\Phi_0((1,0);-;(p,q),z) & = & 
                    \s\f{p^{n(n-1)/2}}{((p,q);(p,q))_n}z^n, \\ 
E_{p,q}(z) = {}_1\Phi_0((0,1);-;(p,q),-z) & = & 
                    \s\f{q^{n(n-1)/2}}{((p,q);(p,q))_n}z^n, 
\label{pqexps}
\eea 
for the $(p,q)$-exponentials such that  
\be
e_{p,q}(z)E_{p,q}(-z) = 1. 
\ee 
For $p$ $=$ $1$, $e_{1,q}(z)$ and $E_{1,q}(z)$ become $e_q(z)$ and 
$E_q(z)$ respectively.  

For $n$ $=$ $3$ the above general result and the relation~(\ref{abba}) 
imply 
\bea 
   &   &  _1\Phi_0((u,v);-;(p,q),z){}_1\Phi_0((v,w);-;(p,q),z) \nn \\ 
   &   &  \qquad \qquad =  {}_1\Phi_0((u,w);-;(p,q),z).
\label{uvw}
\eea
Now, if we take $u=1$, $v=a$, $w=ab$ and $p=1$ then this 
equation~(\ref{uvw}) is just the well known product formula for $_1\phi_0$, 
namely,  
\be
_1\phi_0(a;-;q,z){}_1\phi_0(b;-;q,az) = {}_1\phi_0(ab;-;q,z)\,,
\label{prodrule}
\ee
in view of the relation~(\ref{phi21}).  Thus we get a new way of looking 
at the product formula~(\ref{prodrule}) within the $(p,q)$-series 
formalism.  \\ 

\noindent 
{\bf 4. $(p,q)$-Binomial coefficient} \\

\noindent  
The definition  
\be
\ls\ba{c}
n \\
k 
\ea\rs_{p,q} = \f{((p,q);(p,q))_n}{((p,q);(p,q))_k((p,q);(p,q))_{n-k}},
               \quad k = 0,1,\dots,n,  
\ee
provides a natural generalization of the $q$-binomial coefficient.  In 
terms of the $(p,q)$-number the $(p,q)$-binomial coefficient (written 
without the subscript $p,q$) becomes 
\be
\ls\ba{c}
n \\
k 
\ea\rs = \f{[n]!}{[k]![n-k]!},
\ee
where, as usual, 
\be
[n]! = [n][n-1]\dots [2][1], \qquad [0]! = 1. 
\ee
Then, the result 
\bea
   &   & _1\Phi_0((p^n,q^n);-;(p,q),z) 
         = \sum_{k=0}^\infty\ls\ba{c}
                                n-1+k \\
                                  k 
                               \ea\rs z^k \nn \\
   &   & \qquad = \f{p^{n(n+1)/2}}{((p,p^nz);(p,q))_n} 
         = \lc\sum_{k=0}^n\ls\ba{c}
                               n \\
                               k 
                             \ea\rs (pq)^{k(k-1)/2}(-z)^k\rc^{-1},
\label{pqbin}
\eea 
follows by taking $a = p^n$ and $b = q^n$ in~(\ref{abba}).  The 
relation~(\ref{pqbin}) is obviously a generalization of the result 
\bea 
_1\phi_0(q^n;-;q,z) 
   & = & \sum_{k=0}^\infty\ls\ba{c}
                         n-1+k \\
                           k 
                         \ea\rs_q z^k 
     =  \f{1}{(z;q)_n} \nn \\ 
   & = & \lc\sum_{k=0}^n\ls\ba{c}
                         n \\
                         k 
                       \ea\rs_q q^{k(k-1)/2}(-z)^k\rc^{-1}.
\eea
If we take $p$ $=$ $0$ in~(\ref{pqbin}) we get, 
correctly of course, 
\be
\sum_{k=0}^\infty(q^{n-1}z)^k = \f{1}{1-q^{n-1}z}.
\ee
It should be noted that there is no analogue for the choice $p$ $=$ $0$ 
in the usual $q$-series formalism.  We can also take the limit $p$ $\lra$ 
$q\neq 1$.  Then, the equation~(\ref{pqbin}) takes the form 
\bea
   &   & _1F_0(n;-;q^{n-1}z) 
         = \sum_{k=0}^\infty\lb\ba{c}
                                n-1+k \\
                                  k
                                \ea\rb(q^{n-1}z)^k \nn \\
   &   & \qquad = (1-q^{n-1}z)^{-n} 
                = \lc\sum_{k=0}^n\lb\ba{c}
                                     n \\
                                     k 
                                    \ea\rb(-q^{n-1}z)^k\rc^{-1}.
\eea
Thus, it is seen that, though a $(p,q)$-identity may be derived starting 
with a $q$-identity, the $(p,q)$-identity offers more choices for 
manipulations.  If we choose $(p,q)$ $=$ $(q^{-1},q)$, then, the 
identity~(\ref{pqbin}) becomes 
\bea
   &   & _1\Phi_0((q^{-n},q^n);-;(q^{-1},q),z) = 
          \sum_{k=0}^\infty\ls\ba{c}
                               n-1+k \\
                                 k 
                              \ea\rs_{q^{-1},q}z^k \nn \\
   &   & \qquad = \f{q^{-n(n+1)/2}}{((q^{-1},zq^{-n});(q^{-1},q))_n} 
                = \lc\sum_{k=0}^n\ls\ba{c}
                                      n \\
                                      k 
                                    \ea\rs_{q^{-1},q}(-z)^k\rc^{-1},
\eea
with 
\be
\ls\ba{c}
n \\
k 
\ea\rs_{q^{-1},q} = 
\f{((q^{-1},q);(q^{-1},q))_n}{((q^{-1},q);
    (q^{-1},q))_k((q^{-1},q);(q^{-1},q))_{n-k}}, \quad k = 0,1,\dots,n. 
\ee
which should be relevant in the context of quantum groups. \\

From~(\ref{pqbin}), let us take 
\be
p^{-n(n+1)/2}((p, p^nz);(p,q))_n 
    = \sum_{k=0}^n\ls\ba{c}
                      n \\
                      k 
                     \ea\rs(pq)^{k(k-1)/2}(-z)^k.
\label{pnz}
\ee
Using~(\ref{ll}) and taking $z = \z_q/\z_p$, we can rewrite~(\ref{pnz}) as 
\be
((p\zeta_p,p^n\zeta_q);(p,q))_n 
      = \sum_{k=0}^n\ls\ba{c}
                        n \\
                        k 
                       \ea\rs p^{(n(n+1)+k(k-1))/2}q^{k(k-1)/2}(-1)^k 
                              \z_q^k\z_p^{n-k}.
\ee
Now, renaming $p\z_p$ and $p^n\z_q$ as $a$ and $b$, respectively, we get 
\be
((a,b);(p,q))_n 
        = \sum_{k=0}^n\ls\ba{c}
                          n \\
                          k 
                         \ea\rs(-1)^kp^{(n-k)(n-k-1)/2}q^{k(k-1)/2}
                               a^{n-k}b^k.
\label{gbin}
\ee
The $(p,q)$-binomial theorem derived in~\cite{KK}, using the recursion 
relations of the $(p,q)$-binomial coefficients, corresponds to~(\ref{gbin}) 
with the notations $a$ $=$ $\l$, $b$ $=$ $-x$. 

An operator, or noncommutative, form of the $q$-binomial theorem is 
known~\cite{GR}:  If $x$ and $y$ are noncommuting variables such that 
$xy$ $=$ $qyx$, $q$ commutes with $x$ and $y$, and the associative law 
holds, then 
\be
(x+y)^n = \sum_{k=0}^n\ls\ba{c}
                          n \\
                          k 
                         \ea\rs_qy^kx^{n-k} 
       = \sum_{k=0}^n\ls\ba{c}
                          n \\
                          k 
                         \ea\rs_{q^{-1}}x^ky^{n-k}.
\label{opqbin}
\ee
A $(p,q)$-extension of this result is derived in~\cite{SW}, in a specific 
context of a quantum group.  This result can be stated in a general form 
as follows:
\be
(ax+by)^n = \sum_{k=0}^n\ls\ba{c}
                            n \\
                            k 
                           \ea\rs_{p,q}a^{n-k}b^ky^kx^{n-k} 
         = \sum_{k=0}^n\ls\ba{c}
                            n \\
                            k 
                           \ea\rs_{p^{-1},q^{-1}}b^{n-k}a^kx^ky^{n-k},
\label{oppqbin}
\ee
where $ab$ $=$ $p^{-1}ba$, $xy$ $=$ $qyx$, all other commutators among 
the variables $\{a,b,x,y\}$ vanish, $p$ and $q$ commute with $\{a,b,x,y\}$, 
and the associative law holds.  Proof of~(\ref{oppqbin}) follows by replacing 
in~(\ref{opqbin}) $q$ by $q/p$ and $(x,y)$ by $(ax,by)$, and reexpressing the 
result in terms of $p$, $q$, $a$, $b$, $x$, and $y$.  In deriving the second 
part of~(\ref{oppqbin}) one has to use the formula
\be
((a,b);(p,q))_n = 
         (-1)^na^nb^n(pq)^{n(n-1)/2}((a^{-1},b^{-1});(p^{-1},q^{-1}))_n.
\ee 

\noindent
{\bf 5. $(p,q)$-Heine transformation for $_2\Phi_1$} \\

\noindent  
The Heine transformation of the $_2\phi_1$ series, namely, 
\be
_2\phi_1(a,b;c;q,z) = \f{(b,az;q)_\infty}{(c,z;q)_\infty} 
{}_2\phi_1(c/b,z;az;q,b)\,,
\label{qht}
\ee
has the following $(p,q)$-analogue:
\bea
   &   & _2\Phi_1((a,b),(c,d);(e,f);(p,q),z) \nn \\ 
   &   & \ =  
\f{((ce,de),(pe,bcz);(p,q))_\infty}{(ce,cf),(pe,acz);(p,q))_\infty} \nn \\ 
  &    &  \qquad \times 
          {}_2\Phi_1((de,cf),(pe,acz);(pe,bcz);(p,q),p/ce), 
\label{pqht}
\eea

\noindent{\em Proof}: By the Heine transformation~(\ref{qht})  
\bea
   &   &  _2\phi_1(b/a,d/c;f/e;q/p,\z) \nn \\ 
   &   &  \qquad = \f{(d/c,b\z/a;q/p)_\infty}{(f/e,\z;q/p)_\infty}
                   {}_2\phi_1(cf/de,\z;b\z/a;q/p,d/c).
\eea
Using~(\ref{phi12}) and following arguments of the type used in~(\ref{lim}) 
we can rewrite this equation as
\bea
   &   & _2\Phi_1((a,b),(c,d);(e,f);(p,q),pe\z/ac) \nn \\ 
   &   & \qquad = \f{((c,d),(a,b\z);(p,q))_\infty}
             {((e,f),(ac/e,ac\z/e);(p,q))_\infty}  \nn \\ 
   &   & \qquad {}_2\Phi_1((de,cf),(1,\z);(a,b\z);(p,q),pa/ce). 
\eea
Now, taking $\z$ $=$ $acz/pe$, we get 
\bea
   &   & _2\Phi_1((a,b),(c,d);(e,f);(p,q),z) \nn \\ 
   &   & \qquad  = \f{((ce,de),(pe,bcz);(p,q))_\infty}
                     {(ce,cf),(pe,acz);(p,q))_\infty} \nn \\ 
   &   & \qquad {}_2\Phi_1((de,cf),(pe,acz);(pe,bcz);(p,q),p/ce), 
\eea
thus, arriving at the $(p,q)$-Heine transformation formula~(\ref{pqht}) for 
$_2\Phi_1$. \\ 

Setting $a$ $=$ $0$, $b$ $=$ $c$ $=$ $e$ $=$ $1$, relabeling $d$ as $a$ 
and $f$ as $b$, and taking $p$ $=$ $1$, in~(\ref{pqht}) we obtain the 
transformation 
\be
_1\phi_1(a;b;q,z) = \f{(a,z;q)_\infty}{(b;q)_\infty} 
{}_2\phi_1(0,b/a;z;q,a)\,,
\label{phi11}
\ee
which can be directly derived from the $q$-Heine transformation 
formula~(\ref{qht}) by using the limiting process of confluence, namely, 
replacing $z$ by $z/a$ and taking the limit $a$ $\lra$ $\infty$, and then 
relabeling the parameters.  Now, taking $z$ $=$ $b/a$ in~(\ref{phi11}) 
one obtains, using the $q$-binomial theorem, the summation formula~\cite{GR} 
\be
_1\phi_1(a;b;q,b/a) = \f{(b/a;q)_\infty}{(b;q)_\infty}\,,
\label{phi11id}
\ee
which can also be obtained from the $(p,q)$-Gauss sum~(\ref{pqgauss}), 
given below, with the same choice of parameters. \\

\noindent{\bf 6. $(p,q)$-Gauss sum} \\

\noindent 
Using the Heine transformation~(\ref{qht}) one obtains the $q$-Gauss sum
\be
_2\phi_1(a,b;c;q,c/ab) = \f{(c/a,c/b;q)_\infty}{(c,c/ab;q)_\infty}, 
                         \qquad |c/ab| < 1.
\label{qgauss}
\ee
The $(p,q)$-Gauss sum takes the form 
\bea
   &   & {}_2\Phi_1((a,b),(c,d);(e,f);(p,q),pf/bd) \nn \\ 
   &   & \qquad \qquad = \f{((be,af),(de,cf);(p,q))_\infty}
                           {((e,f),(bde,acf);(p,q))_\infty}, 
         \qquad |acf/bde| < 1.
\label{pqgauss}
\eea

\noindent{\em Proof}: Let $z$ $=$ $pf/bd$ in the $(p,q)$-Heine transformation 
formula~(\ref{pqht}).  The result is 
\bea
   &   & _2\Phi_1((a,b),(c,d);(e,f);(p,q),pf/bd) \nn \\ 
   &   & \qquad =  \f{((ce,de),(pe,pcf/d);(p,q))_\infty}
                   {((ce,cf),(pe,pacf/bd);(p,q))_\infty} \nn \\
   &   & \qquad \qquad \times 
         {}_2\Phi_1((de,cf),(pe,pacf/bd);(pe,pcf/d);(p,q),p/ce), \nn \\
   &   & \qquad = \f{((ce,de),(bde,bcf);(p,q))_\infty}
                      {((ce,cf),(bde,acf);(p,q))_\infty} \nn \\ 
   &   & \qquad \qquad \times 
         {}_2\Phi_1((de,cf),(pe,pacf/bd);(pe,pcf/d);(p,q),p/ce). 
\eea 
Note that 
\bea
   &   & _2\Phi_1((de,cf),(pe,pacf/bd);(pe,pcf/d);(p,q),p/ce) \nn \\  
   &   & \qquad = {}_2\Phi_1((de,cf),(bde,acf);(bde,bcf);(p,q),p/ce) \nn \\
   &   & \qquad = {}_1\Phi_0((bde,acf);-;(p,q),p/bce) \nn \\
   &   & \qquad = \f{((p,paf/be);(p,q))_\infty}{((p,pd/c);(p,q))_\infty} 
                = \f{((be,af);(p,q))_\infty}{((be,bde/c);(p,q))_\infty}, 
\eea
in view of the $(p,q)$-binomial theorem and~(\ref{lim}).  Hence, 
\bea
   &   & _2\Phi_1((a,b),(c,d);(e,f);(p,q),pf/bd) \nn \\ 
   &   &  \qquad = \f{((ce,de),(bde,bcf),(be,af);(p,q))_\infty}
                     {((ce,cf),(bde,acf),(be,bde/c);(p,q))_\infty} \nn \\ 
   &   &  \qquad = \f{((c,d),(de,cf),(be,af);(p,q))_\infty}
                     {((e,f),(bde,acf),(c,d);(p,q))_\infty}  \nn \\ 
   &   &  \qquad = \f{((de,cf),(be,af);(p,q))_\infty}
                     {((e,f),(bde,acf);(p,q))_\infty}.  
\eea
Thus, the $(p,q)$-Gauss sum~(\ref{pqgauss}) is derived. 

The identity 
\be
\s\f{q^{n^2}}{(q,qz;q)_n}z^n = \f{1}{(qz;q)_\infty},
\ee
is usually obtained from the $q$-Gauss sum~(\ref{qgauss}) by setting $c$ 
$=$ $qz$ and letting $a$ $\lra$ $\infty$ and $b$ $\lra$ $\infty$.  It should 
be noted that this identity follows immediately from the  $(p,q)$-Gauss 
sum~(\ref{pqgauss}) by mere substitution $a$ $=$ $c$ $=$ $0$, $b$ $=$ $d$ 
$=$ $e$ $=$ $1$, $f$ $=$ $qz$ and $p$ $=$ $1$. 

Another useful form of~(\ref{pqgauss}) is 
\bea
   &   &  _2\Phi_1((a,1),(b,c);(d,\sg c);(p,q),\sg p) \nn \\ 
   &   & \qquad  = \f{((d,\sg ac),(d,\sg b);(p,q))_\infty}
              {((d,\sg c),(d,\sg ab);(p,q))_\infty}.  
         \qquad |\sg ab/d| < 1. 
\label{sigma}
\eea
Now, substituting in~(\ref{sigma}) $a$ $=$ $c$ $=$ $0$, $b$ $=$ $d$ $=$ $1$, 
$\sg$ $=$ $\sqrt{q} z$ and $p$ $=$ $1$, one gets another well-known identity 
\be
\s\f{(-1)^n q^{n^2/2}}{(q;q)_n}z^n = (\sqrt{q}z;q)_\infty, 
\ee
which is usually obtained from the $q$-Gauss sum~(\ref{qgauss}) by setting 
$c$ $=$ $\sqrt{q} bz$ and then letting $b$ $\lra$ $0$ and $a$ $\lra$ $\infty$.  
These examples illustrate the usefulness of the $(p,q)$-series formalism 
even for the treatment of the usual $q$-series.  \\ 

\noindent 
{\bf 7. $(p,q)$-Ramanujan sum} \\

\noindent
Let us assume the obvious $(p,q)$-generalizations of the basic notations and 
definitions associated with bilateral $q$-hypergeometric series.  Thus, 
we write 
\bea
((a,b);(p,q))_{-n} 
        & = & \f{1}{((ap^{-n},bq^{-n});(p,q))_n} \nn \\
        & = & \f{1}{(ap^{-1}-bq^{-1})(ap^{-2}-bq^{-2})\dots
                               (ap^{-n}-bq^{-n})} \nn \\ 
       & = & \f{(-pq/ab)^n(pq)^{n(n-1)/2}}{((p/a,q/b);(p,q))_n}, 
\eea
and 
\bea
   &   & _1\Psi_1((a,b);(c,d);(p,q),z) 
           = \ds\f{((a,b);(p,q))_n}{((c,d);(p,q))_n}z^n   \nn \\
   &   & \qquad = \s\f{((a,b);(p,q))_n}{((c,d);(p,q))_n}z^n + 
                  \sum_{n=1}^\infty \f{((p/c,q/d);(p,q))_n}
                   {((p/a,q/b);(p,q))_n}\lb\f{cd}{abz}\rb^n. 
\eea  
One can show that 
\be
_1\Psi_1((a,b);(c,d);(p,q),z) = {}_1\psi_1(b/a;d/c;q/p,za/c)
\ee
where $_1\psi_1$ is the usual bilateral $q$-series.  Then, using the 
Ramanujan sum, 
\be
{}_1\psi_1(a;b;q,z) = 
\f{(q,b/a,az,q/az;q)_\infty}{(b,q/a,z,b/az;q)_\infty}, \qquad |b/a|<|z|<1, 
\ee
one can show that the $(p,q)$-analogue of the Ramanujan sum is 
\bea 
   &   &  _1\Psi_1((a,b);(c,d);(p,q),z) \nn \\
   &   &  \qquad = \f{((p,q),(bc,ad),(c,bz),(pbz,qc);(p,a))_\infty}
                   {((c,d),(pb,qa),(c,az),(pbz,pd);(p,a))_\infty}, \nn \\
   &   & \qquad \qquad \qquad |ad/bc|<|z|<1. 
\eea
To obtain the $(p,q)$-analogue of the Jacobi triple product identity 
from this the steps are:
(i) $(a,b)\lra(1/a,1/b)$, $z\lra zb/a$, 
(ii) $d=0$, $(p,q)\lra(p^2,q^2)$, $z\lra zq/p$, 
(iii) $b\lra 0$, and 
(iv) $(p,q)\lra(\sqrt{p},\sqrt{q})$.  
The result is: 
\bea 
   &   &  \ds (-1)^n(q/p)^{n^2/2}(z/ac)^n \nn \\ 
   &   &  \qquad = 
\f{((p,q),(\sqrt{p}ca,\sqrt{q}z),(\sqrt{p}z,\sqrt{q}ca);(p,q))_\infty} 
{((p,0),(\sqrt{p}ca,0),(\sqrt{p}z,0);(p,q))_\infty},
\label{pqjtp}
\eea 
which is same as the well known $q$-result with the replacements 
$q\lra q/p$ and $z\lra z/ac$.  The usual Jacobi triplet product identity 
can also be obtained in a simpler way directly from $_1\Psi_1$ by letting 
$a$ $=$ $d$ $=$ $0$, $b$ $=$ $c$ $=$ $1$, $p$ $=$ $1$ and $z\lra z\sqrt{q}$.  

Taking $ac$ $=$ $1$ in~(\ref{pqjtp}), we can also write the $(p,q)$-analogue 
of the Jacobi triple product, for $q<p$, $|z|<1$, as 
\bea  
   &   &  \ds (-1)^n (q/p)^{n^2/2}z^n \nn \\ 
   &   &  \qquad = \prod_{n=1}^\infty     
\f{(p^n-q^n)(p^{n-1/2}-q^{n-1/2}z)(p^{n-1/2}z-q^{n-1/2})}{p^{3n-1}z}. 
\eea  
The Euler identity follows from the $_1\Psi_1$-sum by taking $a$ $=$ $d$ 
$=$ $0$, $b$ $=$ $c$ $=$ $1$, $(p,q)\lra(1,q^3)$, and $z$ $=$ $q$: 
\be
\ds (-1)^nq^{(3n^2-n)/2} = (q;q)_\infty. 
\ee 

\noindent 
{\bf 8. $(p,q)$-Special functions} \\

\noindent 
Let us now make some brief observations on the $(p,q)$-generalizations of 
the $q$-special functions.  First let us consider an example. It is seen 
that  
\be
\ls\ba{c}
    n \\
    k 
   \ea\rs_{p,q} = \ls\ba{c}
                      n \\
                      n-k 
                     \ea\rs_{p,q} 
                = p^{k(n-k)}\ls\ba{c}
                                n \\
                                k 
                               \ea\rs_{q/p} 
                = p^{k(n-k)}\ls\ba{c}
                                n \\
                                n-k
                               \ea\rs_{q/p}.
\label{pqqp}
\ee
The continuous $q$-Hermite polynomial is given by 
\be
H_n(x|q) = \sum_{k=0}^n\ls\ba{c}
                           n \\
                           k 
                          \ea\rs_qe^{i(n-2k)\t}, \qquad x = \cos\t.  
\ee
We may define a continuous $(p,q)$-Hermite polynomial as 
\be
{\cal H}_n(x|p,q) = \sum_{k=0}^n\ls\ba{c}
                                    n \\
                                    k 
                                   \ea\rs_{p,q}e^{i(n-2k)\t}, 
                    \qquad x = \cos\t. 
\ee
In view of the relation~(\ref{pqqp}) it is found that ${\cal H}_n(x|p,q)$ 
is not just $H_n(x|(q/p))$ with a rescaling of $x$: {\em e.g.}, letting 
$(p,q)\lra(q^\a,q^\b)$ one would get a two-parameter family of generalized 
continuous $q$-Hermite polynomials, say $\lc H_n^{(\a,\b)}(x|q)\rc$ with 
the usual $H_n(x|q)$ identified as $H_n^{(0,1)}(x|q)$.  This is in contrast 
to the case of $_r\Phi_{r-1}$ which can always be identified, as already 
noted (see~(\ref{phi12}) and~(\ref{phi21})), with an $_r\phi_{r-1}$; in 
this sense, $_r\Phi_{r-1}$ may be considered a trivial generalization - 
examples in this category would be the $(p,q)$-generalizations of 
$q$-Krawtchouk polynomials, $q$-Meixner polynomials, $q$-Racah polynomials, 
$q$-Askey-Wilson polynomials, $q$-Jacobi polynomials, $q$-Hahn polynomials, 
$q$-Charlier polynomials, continuous $q$-ultraspherical polynomials, etc...  
However, such generalizations are also of interest from the point of view of 
physical applications: one example of such a situation is the study of the 
Clebsch-Gordon coefficients of the two-parameter quantum algebra 
$U_{p,q}(gl(2))$ - a simple relation between the CG-coefficients of 
$U_{p,q}(gl(2))$ and $U_q(sl(2))$ exists~\cite{JV} which must be due to the 
connection between the CG-coefficients of $U_q(sl(2))$ and $_3\phi_2$ 
(see, {\em e.g.},~\cite{SR})).  $(p,q)$-generalizations of gamma and beta 
functions are straightforward~\cite{BK}.  Besides the continuous $q$-Hermite 
polynomials, there are several examples for which the $(p,q)$-generalization 
is nontrivial: discrete $q$-Hermite polynomials, $q$-Laguerre polynomials, 
$q$-Bessel functions ($J_\nu^{(2)}(x;q)$), etc...   We hope to return to 
these topics elsewhere.  \\

\noindent
{\bf 9. Conclusion} \\

\noindent 
We have shown that it is profitable to study the $(p,q)$-hypergeometric 
series, or the twin-basic hypergeometric series, following naturally from 
the extension of the $q$-number $(1-q^n)/(1-q)$ to the twin-basic number 
$(p^n-q^n)/(p-q)$.  In particular, we have studied the $(p,q)$-analogues 
of the $q$-binomial theorem, $q$-binomial coefficient, Heine transformation 
for $_2\phi_1$, Gauss sum for $_2\phi_1$, and the Ramanujan sum for 
$_1\psi_1$.  Further, we have made some brief observations on the 
$(p,q)$-generalizations of the $q$-special functions.  In general, we 
have noted that many of the $q$-results can be generalized directly to 
$(p,q)$-results and once we have the $(p,q)$-results the $q$-results can 
be obtained more easily by mere substitutions for the parameters instead 
of any limiting process as required in the usual $q$-theory.  We believe 
that a detailed study of the $(p,q)$-hypergeometric series, or the 
twin-basic hypergeometric series, should be very interesting. \\     

\noindent 
\small{{\em Acknowledgements}: 
One of us (R.J) is thankful to the organizing committee of the International 
Conference on Number Theory and Mathematical Physics, held at the Srinivasa 
Ramanujan Centre, SASTRA, on the occasion of Ramanujan birth anniversary,  
for the invitation to present this work in the conference.  R.J wishes also 
to thank SASTRA for the kind hospitality during the conference}.

\end{document}